\documentclass[11pt,a4paper]{article}

\title{\bf A note on Gaussian correlation inequalities for nonsymmetric sets}
\author{Adrian P. C. Lim$^a$,\quad Dejun Luo$^{a,b}$\footnote{Email: luodj@amss.ac.cn}
\vspace{3mm}\\
{\footnotesize $^a$UR Math\'{e}matiques, Universit\'{e} de
Luxembourg, 6, rue Richard Coudenhove-Kalergi, L-1359 Luxembourg}\\
{\footnotesize $^b$Key Lab of Random Complex Structures and Data
Science, Academy of Mathematics and}\\
{\footnotesize  Systems Science, Chinese Academy of Sciences,
Beijing 100190, China} }
\date{}

\usepackage{amssymb,amsmath,amsfonts,amsthm,color}

\def\R{\mathbb{R}}
\def\C{\mathcal{C}}
\def\d{\textup{d}}
\def\supp{\textup{supp}}
\def\dist{\textup{dist}}
\def\da{\downarrow}
\def\ra{\rightarrow}
\def\<{\langle}
\def\>{\rangle}

\def\fin{\hfill$\square$}

\newtheorem{theorem}{Theorem}[section]
\newtheorem{lemma}[theorem]{Lemma}
\newtheorem{corollary}[theorem]{Corollary}

\newtheorem{remark}[theorem]{Remark}

\begin{document}

\maketitle
\makeatletter % '@' is now a normal "letter" for TeX
\renewcommand\theequation{\thesection.\arabic{equation}}
\@addtoreset{equation}{section}
\makeatother % '@' is restored as a "non-letter" character for TeX

\begin{abstract}
We consider the Gaussian correlation inequality for nonsymmetric
convex sets. More precisely, if $A\subset\R^d$ is convex and the
origin $0\in A$, then for any ball $B$ centered at the origin, it
holds $\gamma_d(A\cap B)\geq \gamma_d(A)\gamma_d(B)$, where
$\gamma_d$ is the standard Gaussian measure on $\R^d$. This
generalizes Proposition 1 in [Arch. Rational Mech. Anal. 161 (2002),
257--269].
\end{abstract}



\section{Introduction}

Let $A,B\subset\R^d$ be symmetric convex subsets. The Gaussian
correlation inequality claims that
  \begin{equation}\label{correl}
  \gamma_d(A\cap B)\geq \gamma_d(A)\gamma_d(B),
  \end{equation}
where $\gamma_d$ is the standard Gaussian measure on $\R^d$. The
case $d=1$ is trivial, since $A$ and $B$ are centered intervals,
hence one is contained in the other. The case $d=2$ was proved by
Pitt \cite{Pitt77}. For higher dimensional cases, there are only
partial results. For instance, in \cite{SSZ} it was shown that
\eqref{correl} holds if $A$ and $B$ are ellipsoids, which was soon
generalized by Harg\'e \cite{Harge} to allow one of them to be an
arbitrary symmetric convex set. Harg\'e's proof relies on the
modified Ornstein-Uhlenbeck semigroup and the properties of
log-concave functions. A rather short proof of Harg\'e's result was
presented in \cite{Cordero}, by making use of the deep results in
the theory of optimal transport, that is, the optimal transport map,
which pushes forward the Gaussian measure $\gamma_d$ to a
probability measure $\nu$ having a log-concave density with respect
to $\gamma_d$, is a contraction (see \cite{Caffarelli}). In
\cite{Hu97}, the author obtained a correlation inequality for the
Gaussian measure via a formula for It\^o-Wiener chaos expansion. Li
W.V. \cite{Li99} presented a weaker form of the correlation
inequality \eqref{correl}, which is useful to show the existence of
small ball constants. For a more detailed survey of the studies on
\eqref{correl}, see \cite[Section 2.4]{LiShao}.

In this note we consider two special cases of the correlation
inequality. It is clear that we only need to consider bounded
subsets of $\R^d$. In the sequel, we always assume that the sets are
bounded and closed. First we prove the following result.

\begin{theorem}\label{sect-1-thm.1}
Let $\d\mu=\rho(|x|)\,\d x$ be a probability measure on $\R^d$ with
$\rho\in C(\R_+,(0,\infty))$. Suppose $A\subset\R^d$ is convex and
the origin $0\in A$. Then for any ball $B$ centered at the origin,
we have
  $$ \mu(A\cap B)\geq \mu(A)\mu(B).$$
\end{theorem}

Clearly the Gaussian measure $\gamma_d$ is a special case of $\mu$
considered above. The new point here is that the set $A$ does not
have to be symmetric, at the price of the regularity on $B$. Theorem
\ref{sect-1-thm.1} is a slight generalization of \cite[Proposition
1]{Cordero}; the latter requires that the origin $0$ is the unique
fixed point of all the isometries which leave $A$ (globally)
invariant. Remark also that our proof (see Section 2) uses purely
elementary analysis, while the one in \cite{Cordero} relies on the
result in the theory of optimal transport.

If we want to prove the correlation inequality for more general sets
$B$ other than the balls (e.g. the ellipsoids), then some additional
conditions have to be imposed on the set $A$.

\begin{theorem}\label{sect-1-thm.2}
Assume that $\mu$ is a product probability measure:
$\mu=\prod_{i=1}^d\mu_i$, where $\d\mu_i=\rho_i(|x_i|)\,\d x_i$ with
$\rho_i\in C(\R_+,(0,\infty))$. Let $A\subset \R^d$ be a convex set
with the following property: $x\in A$ implies that its projections
on all the coordinate hyperplanes also belong to $A$. Then for any
ellipsoid
  $$B=\bigg\{x\in\R^d:\frac{x_1^2}{a_1^2}+\cdots+\frac{x_d^2}{a_d^2}\leq 1\bigg\},$$
where $a_1,\cdots,a_d$ are positive constants, we have $\mu(A\cap
B)\geq \mu(A)\,\mu(B)$.
\end{theorem}

This result will be proved in Section 3. It is easy to see that the
set $A$ considered in Theorem \ref{sect-1-thm.2} contains the origin
$0$. An example for the set $A$ is $\{x=(x_1,\cdots,x_d)\in\R^d:
\forall\,i=1,\cdots,d,\, x_i\geq0 \mbox{ and } \sum_{i=1}^dx_i\leq
1\}$. We would like to mention that Theorem \ref{sect-1-thm.2} still
holds for more general set $B$, see Remark \ref{sect-3-rem-1}.

A nonnegative function $f:\R^d\ra\R_+$ is called log-concave if for
any $x,y\in\R^d$ and $0<\lambda<1$, $f(\lambda x+(1-\lambda)y)\geq
f(x)^\lambda f(y)^{1-\lambda}$. For any convex set $A\subset \R^d$,
one easily concludes that the indicator function ${\bf 1}_A$ is
log-concave. The Gaussian correlation inequality \eqref{correl} has
the following functional version: for any log-concave and symmetric
functions $f,g$, it holds
  \begin{equation}\label{correl.1}
  \gamma_d(fg)\geq \gamma_d(f) \gamma_d(g).
  \end{equation}
Here $\gamma_d(f)=\int_{\R^d}f\,\d\gamma_d$. Following the method of
\cite[Section 3]{Cordero}, we show in the last section that
\eqref{correl.1} holds if $f$ is log-concave and $g=\varphi(\<
\Sigma x,x\>)$, where $\varphi\in C(\R_+,\R_+)$ is a decreasing
function and $\Sigma$ is a positive definite matrix. So we give an
alternative proof to \cite[Theorem 2]{Harge}. By an approximation
argument, we obtain again Harg\'e's result.

\section{Proof of Theorem \ref{sect-1-thm.1}}

In this section we prove Theorem \ref{sect-1-thm.1}. In fact we will
prove a more general result. To this end, we introduce a class of
functions on $\R^d$:
  \begin{equation}\label{funct-class}
  \C_d=\big\{f\in C_c(\R^d,\R_+):\forall\, c>0, \{f>c\} \mbox{ is convex and }
  \forall\,x\in\R^d, f(x)\leq f(0)\big\}.
  \end{equation}
Let $S^{d-1}$ be the unit sphere in $\R^d$. For a bounded measurable
function $g$ on $\R^d$, define $\mu(g)=\int_{\R^d}g\,\d\mu$. We have
the following simple observations.

\begin{lemma}\label{sect-2-lem.1}
Let $f\in\C_d$ and $f\neq0$. Then
\begin{itemize}
\item[\rm(i)] for any $\theta\in S^{d-1}$, the function $t\mapsto f(t\theta)$
is decreasing on $\R_+$;
\item[\rm(ii)] $f(0)>\mu(f)$.
\end{itemize}
\end{lemma}

\noindent{\bf Proof.} (i) Suppose that there are $t_1<t_2$ such that
$f(t_1\theta)<f(t_2\theta)$. Consider the set
$E:=\big\{f>[f(t_1\theta)+ f(t_2\theta)]/2\big\}$. Then
$0,t_2\theta\in E$ but $t_1\theta\in E^c$, which contradicts the
fact that $E$ is convex.

(ii) By the definition of the class $\C_d$,
  $$f(0)-\mu(f)=\int_{\R^d}(f(0)-f(x))\,\d\mu(x)\geq0.$$
If $f(0)=\mu(f)$, then $f(x)\equiv f(0)$ for all $x\in\R^d$, which
is impossible. Hence $f(0)>\mu(f)$. \fin

\medskip

Now we prove

\begin{theorem}\label{sect-2-thm.1}
Assume that $\d\mu=\rho(|x|)\,\d x$ is a probability measure on
$\R^d$ with $\rho\in C(\R_+,(0,\infty))$. For any $f\in \C_d$ and
any ball $B$ centered at the origin, it holds
  $$\mu(f{\bf 1}_B)\geq \mu(f)\,\mu(B).$$
\end{theorem}

\noindent{\bf Proof.} Obviously we can assume $\mu(f)>0$. For
$t\geq0$, let $B_t$ be the ball centered at the origin with radius
$t$. Define the function
  $$\Phi(t)=\mu(f{\bf 1}_{B_t})-\mu(f)\mu(B_t),\quad t\geq0.$$
First we show that $\Phi$ is positive when $t$ is sufficiently small
and large. By Lemma \ref{sect-2-lem.1}(ii), there is $t_0>0$ such
that for all $x\in B_{t_0}$, $f(x)>\mu(f)$. Thus for any
$t\in(0,t_0)$,
  $$\Phi(t)=\int_{B_t}\big[f(x)-\mu(f)\big]\,\d\mu(x)>0.$$
When $t$ is big enough such that $\supp(f)\subset B_t$, we have
  $$\Phi(t)=\mu(f)-\mu(f)\mu(B_t)=\mu(f)\mu(B_t^c)>0.$$

Next we compute the derivative $\Phi^\prime(t)$. We have for $h>0$,
  $$\mu(f{\bf 1}_{B_{t+h}})-\mu(f{\bf 1}_{B_t})=
  \int_{B_{t+h}\setminus B_t}f(x)\,\d\mu(x)
  =\int_{B_{t+h}\setminus B_t}f(x)\rho(|x|)\,\d x.$$
Using the spherical coordinate, the above equality can be written as
  \begin{align*}
  \mu(f{\bf 1}_{B_{t+h}})-\mu(f{\bf 1}_{B_t})&=
  \int_t^{t+h}\bigg(\int_{S^{d-1}}f(r\theta)\rho(|r\theta|)\,\d\sigma(\theta)\bigg)r^{d-1}\,\d
  r\cr
  &=\int_t^{t+h}\bigg(\int_{S^{d-1}}f(r\theta)\,\d\sigma(\theta)\bigg)\rho(r)r^{d-1}\,\d
  r,
  \end{align*}
where $\sigma$ is the volume measure on $S^{d-1}$. Since the
functions $f$ and $\rho$ are continuous, dividing both sides by $h$
and letting $h\ra0$ lead to
  $$\frac{\d}{\d t}\mu(f{\bf 1}_{B_t})=\rho(t)t^{d-1}\int_{S^{d-1}}f(t\theta)\,\d\sigma(\theta).$$
Similarly we have $\frac{\d}{\d t}\mu(B_t)=\rho(t)t^{d-1}
\sigma(S^{d-1})$. Therefore
  \begin{align}\label{sect-2-thm.1.1}
  \Phi^\prime(t)&=\rho(t)t^{d-1}\int_{S^{d-1}}f(t\theta)\,\d\sigma(\theta)
  -\mu(f)\rho(t)t^{d-1}\sigma(S^{d-1})\cr
  &=\rho(t)t^{d-1}\int_{S^{d-1}}\big[f(t\theta)-\mu(f)\big]\,\d\sigma(\theta).
  \end{align}
From this expression, it is clear that $\Phi^\prime$ is continuous.
For $t>0$ small enough, we conclude from \eqref{sect-2-thm.1.1} and
Lemma \ref{sect-2-lem.1}(ii) that $\Phi^\prime(t)>0$. Let
$t_1=\inf\{t>0:\Phi^\prime(t)=0\}$. Then
$\int_{S^{d-1}}\big[f(t_1\theta)-\mu(f)\big]\,\d\sigma(\theta)=0$.
By Lemma \ref{sect-2-lem.1}(i), for any $t>t_1$,
  $$\int_{S^{d-1}}\big[f(t\theta)-\mu(f)\big]\,\d\sigma(\theta)
  \leq \int_{S^{d-1}}\big[f(t_1\theta)-\mu(f)\big]\,\d\sigma(\theta)=0.$$
Hence $\Phi^\prime(t)\leq 0$. It follows that $\Phi(t)$ is
increasing on $[0,t_1]$ and decreasing on $[t_1,\infty)$. Combining
this with the fact that $\Phi(t)>0$ when $t$ is sufficiently small
and large, we complete the proof. \fin

\medskip

Now we are ready to prove Theorem \ref{sect-1-thm.1}.

\medskip

\noindent{\bf Proof of Theorem \ref{sect-1-thm.1}.} We will
construct a sequence of functions, belonging to $\C_d$, which
converge to ${\bf 1}_A$. Let $\dist(\cdot,A)$ be the distance
function to $A$. For $n\geq 1$, define
  \begin{equation}\label{sect-2-proof.1}
  f_n(x)=1-n\big[n^{-1}\wedge\dist(x,A)\big],\quad x\in\R^d.
  \end{equation}
Then it is clear that $f_n\in C_c(\R^d,\R_+)$, $0\leq f_n\leq 1$ and
the restriction $f_n|_A\equiv 1$. Since $0\in A$, we have for all
$x\in \R^d$, $f_n(x)\leq 1=f_n(0)$. It remains to show that for any
$c\in[0,1)$, $\{f_n>c\}$ is convex. In fact,
  $$\{f_n>c\}=\{x\in\R^d:\dist(x,A)<(1-c)/n\}.$$
If $x,y\in \{f_n>c\}$, then $\dist(x,A)\vee \dist(y,A)<(1-c)/n$.
Thus there are $x_0,y_0\in A$ such that $|x-x_0|\vee |y-y_0|<
(1-c)/n$. For any $\lambda\in(0,1)$, we have by the convexity of
$A$, $\lambda x_0+(1-\lambda)y_0\in A$. Moreover,
  $$\big|(\lambda x+(1-\lambda)y)-(\lambda x_0+(1-\lambda)y_0)\big|
  \leq \lambda|x-x_0|+(1-\lambda) |y-y_0|<(1-c)/n.$$
Therefore $\dist\big(\lambda x+(1-\lambda)y,A\big)<(1-c)/n$;
equivalently, $\lambda x+(1-\lambda)y\in \{f_n>c\}$. This means that
$\{f_n>c\}$ is convex.

Now applying Theorem \ref{sect-2-thm.1} to $f_n$, we have
  $$\mu(f_n{\bf 1}_B)\geq \mu(f_n)\,\mu(B),\quad \mbox{for all }n\geq 1.$$
Since $f_n\da{\bf 1}_A$ on $\R^d$, by the monotone convergence
theorem, letting $n\ra\infty$ completes the proof. \fin

\section{Proof of Theorem \ref{sect-1-thm.2}}

In order to prove Theorem \ref{sect-1-thm.2}, we introduce another
family of functions:
  \begin{align*}
  \bar\C_d=\big\{f\in C_c(\R^d,\R_+):&\, \forall\,i\in\{1,\cdots,d\} \mbox{ and }
  (x_1,\cdots,x_{i-1},x_{i+1},\cdots,x_d)\in\R^{d-1} \mbox{ fixed,}\cr
  &\, \mbox{the function } x_i\mapsto f(x_1,\cdots,x_{i-1},
  x_i,x_{i+1},\cdots,x_d)\in \C_1\big\},
  \end{align*}
where $\C_1$ is the class of functions defined in
\eqref{funct-class} for $d=1$. We have

\begin{theorem}\label{sect-3-thm.1} Assume that $\mu$ is a product
probability measure: $\mu=\prod_{i=1}^d\mu_i$, where
$\d\mu_i=\rho_i(|x_i|)\,\d x_i$ with $\rho_i\in C(\R_+,(0,\infty))$.
Let $B$ be an ellipsoid:
  $$B=\bigg\{x\in\R^d:\frac{x_1^2}{a_1^2}+\cdots+\frac{x_d^2}{a_d^2}\leq 1\bigg\},$$
where $a_1,\cdots,a_d$ are positive constants. Then for any $f\in
\bar\C_d$, the following inequality holds:
  $$\mu(f{\bf 1}_B)\geq \mu(f)\,\mu(B).$$
\end{theorem}

\noindent{\bf Proof.} We will prove this result by induction on the
dimension $d$. When $d=1$, this theorem is a special case of Theorem
\ref{sect-1-thm.1}. Now suppose that the assertion is true in the
$d-1$ dimensional case. Denote by $\mu_{(d-1)}=\prod_{i=1}^{d-1}
\mu_i$ the product measure on $\R^{d-1}$. By Fubini's theorem,
  \begin{align}\label{sect-3-thm.1.1}
  \mu(f{\bf 1}_B)&=\int_B f\,\d(\mu_{(d-1)}\times\mu_d)\cr
  &=\int_{-a_d}^{a_d}\d\mu_d(x_d)\int_{B_{d-1}(x_d)}
  f(x_1,\cdots,x_{d-1},x_d)\,\d\mu_{(d-1)}(x_1,\cdots,x_{d-1}),
  \end{align}
where $B_{d-1}(x_d)$ is a $d-1$ dimensional ellipsoid:
  $$B_{d-1}(x_d)=\bigg\{(x_1,\cdots,x_{d-1})\in\R^{d-1}:\frac{x_1^2}{a_1^2}+\cdots
  +\frac{x_{d-1}^2}{a_{d-1}^2}\leq 1-\frac{x_d^2}{a_d^2}\bigg\}.$$
Notice that for fixed $x_d\in[-a_d,a_d]$, $f(\cdot,x_d)\in\bar
\C_{d-1}$. Using the induction hypothesis, we have
  $$\int_{B_{d-1}(x_d)}f(x_1,\cdots,x_{d-1},x_d)\,\d\mu_{(d-1)}(x_1,\cdots,x_{d-1})
  \geq \mu_{(d-1)}(f(\cdot,x_d))\,\mu_{(d-1)}(B_{d-1}(x_d)).$$
Therefore by \eqref{sect-3-thm.1.1},
  \begin{align*}
  \mu(f{\bf 1}_B)\geq\int_{-a_d}^{a_d}\mu_{(d-1)}(f(\cdot,x_d))
  \,\mu_{(d-1)}(B_{d-1}(x_d))\,\d\mu_d(x_d).
  \end{align*}
The function $[-a_d,a_d]\ni x_d\mapsto \mu_{(d-1)}(B_{d-1}(x_d))$ is
even and $\mu_{(d-1)}(B_{d-1}(a_d))=0$. We extend it to a function
on $\R$ by setting $\mu_{(d-1)}(B_{d-1}(x_d))\equiv 0$ for
$|x_d|>a_d$. Then the above inequality becomes
  \begin{align}\label{sect-3-thm.1.2}
  \mu(f{\bf 1}_B)&\geq\int_{-\infty}^\infty\mu_{(d-1)}(f(\cdot,x_d))
  \,\mu_{(d-1)}(B_{d-1}(x_d))\,\d\mu_d(x_d)\cr
  &=\bigg(\int_{-\infty}^0+\int_0^\infty\bigg)\mu_{(d-1)}(f(\cdot,x_d))
  \,\mu_{(d-1)}(B_{d-1}(x_d))\,\d\mu_d(x_d).
  \end{align}
We denote by $I_1$ and $I_2$ the two integrals on the right hand
side of \eqref{sect-3-thm.1.2}. Note that the even function
$x_d\mapsto \mu_{(d-1)}(B_{d-1}(x_d))$ is decreasing on $\R_+$. On
the other hand, for any fixed $(x_1,\cdots,x_{d-1})\in \R^{d-1}$, by
the definition of the class $\bar\C_d$ and Lemma
\ref{sect-2-lem.1}(i), the function $x_d\mapsto
f(x_1,\cdots,x_{d-1},x_d)$ is decreasing (resp. increasing) on
$\R_+$ (resp. $\R_-=(-\infty,0]$). Hence the same is true for
$x_d\mapsto \mu_{(d-1)}(f(\cdot,x_d))$. Applying the FKG inequality
(see Lemma \ref{sect-3-lem.2} below) to $2\mu_d$ on $(-\infty,0]$
leads to
  \begin{align}\label{sect-3-thm.1.3}
  I_1&\geq\frac12\bigg(2\int_{-\infty}^0\mu_{(d-1)}(f(\cdot,x_d))
  \,\d\mu_d(x_d)\bigg)
  \bigg(2\int_{-\infty}^0\mu_{(d-1)}(B_{d-1}(x_d))
  \,\d\mu_d(x_d)\bigg)\cr
  &=\bigg(\int_{-\infty}^0\mu_{(d-1)}(f(\cdot,x_d))
  \,\d\mu_d(x_d)\bigg)\bigg(\int_{-\infty}^\infty\mu_{(d-1)}(B_{d-1}(x_d))
  \,\d\mu_d(x_d)\bigg),
  \end{align}
where the equality follows from the symmetry of the integrand the
the measure $\mu_d$. Similarly we have
  $$I_2\geq \bigg(\int_0^\infty\mu_{(d-1)}(f(\cdot,x_d))
  \,\d\mu_d(x_d)\bigg)\bigg(\int_{-\infty}^\infty\mu_{(d-1)}(B_{d-1}(x_d))
  \,\d\mu_d(x_d)\bigg).$$
Combining this with \eqref{sect-3-thm.1.2} and
\eqref{sect-3-thm.1.3}, we conclude that
  \begin{align*}
  \mu(f{\bf 1}_B)&\geq \bigg(\int_{-\infty}^\infty\mu_{(d-1)}(f(\cdot,x_d))
  \,\d\mu_d(x_d)\bigg)\bigg(\int_{-\infty}^\infty\mu_{(d-1)}(B_{d-1}(x_d))
  \,\d\mu_d(x_d)\bigg)\cr
  &=\mu(f)\,\mu(B).
  \end{align*}
Therefore the result holds as well in the $d$ dimensional case. The
proof is complete. \fin

\begin{lemma}[FKG inequality]\label{sect-3-lem.2}
Let $-\infty\leq a,b\leq\infty$ and $\nu$ be a probability measure
on $[a,b]$. Assume $f$ and $g$ are two bounded increasing (or
decreasing) functions on $[a,b]$, then
  $$\int_a^b fg\,\d\nu\geq \bigg(\int_a^b f\,\d\nu\bigg)\bigg(\int_a^b g\,\d\nu\bigg).$$
\end{lemma}

\noindent{\bf Proof.} For any $x,y\in[a,b]$, since both $f$ and $g$
are increasing (or decreasing) functions on $[a,b]$, we have
  $$(f(x)-f(y))(g(x)-g(y))\geq0.$$
As the two functions are bounded, we can integrate the above
inequality on $[a,b]^2$ with respect to $\nu\times\nu$ and obtain
  $$\int_{[a,b]^2}(f(x)-f(y))(g(x)-g(y))\,\d(\nu\times\nu)(x,y)\geq0.$$
Expanding the product gives the desired result. \fin

\begin{remark}\label{sect-3-rem-1}
The proof of Theorem \ref{sect-3-thm.1} works for more general set
$B$. Indeed, for $i=1,\cdots,d$, let $f_i\in C(\R_+,\R_+)$ be a
strictly increasing function such that $f_i(0)=0$. Then the result
of Theorem \ref{sect-3-thm.1} still holds for the set
  $$B=\{x\in\R^d:f_1(|x_1|)+\cdots+f_d(|x_d|)\leq 1\}.$$
Notice that $B$ can even be non-convex. For example, when $d=2$ and
$f_1(t)= f_2(t)=\sqrt t$ for $t\geq0$, then $B=\{x\in\R^2:
\sqrt{|x_1|}+\sqrt{|x_2|}\leq 1\}$ is clearly not convex.
\end{remark}

\medskip

Now we are in the position to prove Theorem \ref{sect-1-thm.2}. We
focus on the case $d\geq 2$ (the case $d=1$ has been proved in
Theorem \ref{sect-1-thm.1}).

\medskip

\noindent{\bf Proof of Theorem \ref{sect-1-thm.2}.} Consider the
approximations $f_n$ of the indicator function ${\bf 1}_A$ defined
in \eqref{sect-2-proof.1}. In order to apply Theorem
\ref{sect-3-thm.1}, we have to show that for every $n\geq 1$,
$f_n\in\bar\C_d$. For simplicity of notations, we assume $i=1$, that
is, for any $x^\prime=(x_2,\cdots,x_d)\in\R^{d-1}$ fixed, we need to
prove that the function $x_1\mapsto f_n(x_1,x^\prime)\in \C_1$,
where $\C_1$ is defined in \eqref{funct-class}. For any $c>0$,
  $$I:=\{x_1\in\R:f_n(x_1,x^\prime)>c\}=\big\{x_1\in\R:
  \dist\big((x_1,x^\prime),A\big)<(1-c)/n\big\}.$$
If $x_1,\bar x_1\in I$, $x_1<\bar x_1$, then
$\dist\big((x_1,x^\prime),A\big) \vee \dist\big((\bar
x_1,x^\prime),A\big)<(1-c)/n$. Since $A$ is convex, as in the proof
of Theorem \ref{sect-1-thm.1}, we can show that
  $$\dist\big(\lambda(x_1,x^\prime)+(1-\lambda)(\bar
  x_1,x^\prime),A\big)<(1-c)/n$$
for all $\lambda\in(0,1)$. That is, $\dist\big((\lambda
x_1+(1-\lambda)\bar x_1,x^\prime), A\big)<(1-c)/n$. Consequently,
for all $\lambda\in(0,1)$, $\lambda x_1+(1-\lambda)\bar x_1\in I$.
This means that $I$ is an interval, hence it is convex.

Next we prove that $0\in I$ whenever $I$ is nonempty. Indeed, if
$x_1\in I$, then $\dist\big((x_1,x^\prime),A\big)<(1-c)/n$. Hence
there is $y=(y_1,y^\prime)\in A$ such that
$|(x_1,x^\prime)-(y_1,y^\prime)|<(1-c)/n$. By the property of $A$,
we have $(0,y^\prime)\in A$. Moreover
  $$|(0,x^\prime)-(0,y^\prime)|\leq
  |(x_1,x^\prime)-(y_1,y^\prime)|<(1-c)/n.$$
Therefore $\dist\big((0,x^\prime),A\big)<(1-c)/n$, that is, $0\in
I$. Now if there is $x_1\in \R$ such that $f_n(x_1,x^\prime)>f_n(0,
x^\prime)$, then consider the interval
  $$\tilde I=\big\{f_n(\cdot,x^\prime)>\big(f_n(x_1,x^\prime)+
  f_n(0,x^\prime)\big)/2\big\}.$$
We have $x_1\in \tilde I$ but $0\in \tilde I^c$, which is a
contradiction with the result that we have just proved. Hence
$f_n(0,x^\prime)\geq f_n(x_1,x^\prime)$ for all $x_1\in\R$.
Therefore the function $x_1\mapsto f_n(x_1,x^\prime)\in \C_1$.
Summing up these arguments, we conclude that $f_n\in \bar\C_d$.

Now applying Theorem \ref{sect-3-thm.1} to $f_n$, we obtain
$\mu(f_n{\bf 1}_B)\geq \mu(f_n)\,\mu(B)$ for all $n\geq1$. Letting
$n\ra\infty$ gives the desired result. \fin

\section{A special case of \eqref{correl.1}}

In the present section, we follow the method in \cite{Cordero} (see
p.265) and prove that the inequality \eqref{correl.1} holds if $g$
is the composition of a decreasing function and a positive definite
quadratic form. This gives an alternative proof to \cite[Theorem
2]{Harge}.

\begin{theorem}\label{sect-4-thm.1}
Assume that $f\in C(\R^d,\R_+)$ is a log-concave and symmetric
function. Let $\Sigma$ be a positive definite matrix and $\varphi\in
C(\R_+, \R_+)$ a decreasing function. Then
  \begin{equation}\label{sect-4.1}
  \int_{\R^d}f(x)\,\varphi(\<\Sigma x,x\>)\,\d\gamma_d(x)\geq
  \bigg(\int_{\R^d}f(x)\,\d\gamma_d(x)\bigg)
  \bigg(\int_{\R^d}\varphi(\<\Sigma x,x\>)\,\d\gamma_d(x)\bigg).
  \end{equation}
\end{theorem}

\noindent{\bf Proof.} Consider the Gaussian probability measure
  $$\d\mu=\frac{1}{(2\pi)^{n/2}\sqrt{\det(\Sigma)}}e^{-\<\Sigma^{-1}x,x\>/2}\d x,$$
where $\det(\Sigma)$ is the determinant of $\Sigma$. Then
\eqref{sect-4.1} is equivalent to
  \begin{equation}\label{sect-4.2}
  \int_{\R^d}f\big(\sqrt{\Sigma^{-1}}\,x\big)\varphi(|x|^2)\,\d\mu(x)
  \geq\bigg(\int_{\R^d}f\big(\sqrt{\Sigma^{-1}}\,x\big)\d\mu(x)\bigg)
  \bigg(\int_{\R^d}\varphi(|x|^2)\,\d\mu(x)\bigg).
  \end{equation}
Since $f$ is log-concave and symmetric, it is easy to see that
$f(x)\leq f(0)$ for all $x\in\R^d$, hence
$C_f:=\int_{\R^d}f\big(\sqrt{\Sigma^{-1}}\,x\big)\d\mu(x)< +\infty$.
We introduce the probability measure $\mu_f$ defined by
  $$\d\mu_f=\frac1{C_f}f\big(\sqrt{\Sigma^{-1}}\,x\big)\d\mu(x).$$
Hence it is sufficient to prove that
  \begin{equation}\label{sect-4.3}
  \int_{\R^d}\varphi(|x|^2)\,\d\mu_f(x)
  \geq\int_{\R^d}\varphi(|x|^2)\,\d\mu(x).
  \end{equation}
Now let $T$ be the optimal transport map pushing forward the
Gaussian measure $\mu$ to $\mu_f$, i.e. $\mu_f=\mu\circ T^{-1}$.
Since the density function $\frac1{C_f}f\big(\sqrt{\Sigma^{-1}}
\,x\big)$ is also log-concave, we deduce from Caffarelli's result
(see \cite{Caffarelli}) that $T$ is a contraction. Moreover by the
symmetry of the density function, we have $T(-x)=-T(x)$;
particularly, $T(0)=0$. Therefore, $|T(x)|\leq |x|$ for all
$x\in\R^d$. As a result,
  \begin{align*}
  \int_{\R^d}\varphi(|x|^2)\,\d\mu_f(x)&=\int_{\R^d}\varphi(|x|^2)\,\d(\mu\circ
  T^{-1})(x)\cr
  &=\int_{\R^d}\varphi(|T(x)|^2)\,\d\mu(x)\geq
  \int_{\R^d}\varphi(|x|^2)\,\d\mu(x),
  \end{align*}
where the last inequality follows from the fact that $\varphi$ is a
decreasing function. \eqref{sect-4.3} is proved. \fin

\medskip

If we take $\varphi(t)=e^{-t/2}$ for $t\geq0$, then Theorem
\ref{sect-4-thm.1} reduces to \cite[Proposition 2]{SSZ} (see p.352).
Moreover by approximating the indicator function of a symmetric
convex set, we can reprove Harg\'e's result \cite{Harge}.

\begin{corollary} Let $A\subset \R^d$ be any symmetric convex set and
$B$ be the ellipsoid $\{x\in\R^d:\<\Sigma x,x\>\leq 1\}$. Then
  $$\gamma_d(A\cap B)\geq\gamma_d(A)\gamma_d(B).$$
\end{corollary}

\noindent{\bf Proof.} We consider again the sequence of
approximating functions $f_n$ defined in \eqref{sect-2-proof.1}.
First we show that $f_n$ is log-concave for every $n\geq1$. Since
the set $A$ is convex, it is easy to see that the distance function
$\dist(\cdot,A)$ is also convex. Hence for any $x,y\in\R^d$ and
$0<\lambda<1$,
  \begin{equation}\label{sect-4.4}
  \dist\big(\lambda x+(1-\lambda)y,A\big)\leq \lambda\,\dist(x,A)
  +(1-\lambda)\,\dist(y,A).
  \end{equation}
In order to show that $f_n(\lambda x+(1-\lambda)y) \geq
f_n(x)^\lambda f_n(y)^{1-\lambda}$, it is enough to consider the
case $f_n(x)\wedge f_n(y)>0$, that is, $\dist(x,A)\vee \dist(y,A)<
1/n$. Therefore by \eqref{sect-4.4},
  $$1-n\,\dist\big(\lambda x+(1-\lambda)y,A\big)
  \geq \lambda\,[1-n\,\dist(x,A)]+(1-\lambda)\,[1-n\,\dist(y,A)].$$
In other words,
  \begin{equation}\label{sect-4.5}
  f_n(\lambda x+(1-\lambda)y)\geq \lambda\,f_n(x)+(1-\lambda)\,f_n(y).
  \end{equation}
Now using Young's inequality (for any $a,b\geq0$ and $p,q>1$ such
that $p^{-1}+q^{-1}=1$, it holds $ab\leq \frac{a^p}p+\frac{b^q}q$),
we obtain
  $$\lambda\,f_n(x)+(1-\lambda)\,f_n(y)\geq f_n(x)^\lambda f_n(y)^{1-\lambda}.$$
Combining this with \eqref{sect-4.5}, we obtain the log-concavity of
$f_n$.

The symmetry of the set $A$ implies that $f_n$ is also symmetric.
Now applying Theorem \ref{sect-4-thm.1} to the functions $f_n$ and
letting $n\ra\infty$, we arrive at
  \begin{equation}\label{sect-4.6}
  \int_{\R^d}{\bf 1}_A(x)\,\varphi(\<\Sigma x,x\>)\,\d\gamma_d(x)\geq
  \gamma_d(A)\bigg(\int_{\R^d}\varphi(\<\Sigma x,x\>)\,\d\gamma_d(x)\bigg).
  \end{equation}
Next define
  $$\varphi_n(t)=\begin{cases}
  1,& t\in[0,1];\\
  1-n(t-1),& 1<t<1+n^{-1};\\
  0,& t\geq 1+n^{-1}.
  \end{cases}$$
Then $\varphi_n(t)\da {\bf 1}_{[0,1]}(t)$ as $n\ra\infty$. Replacing
$\varphi$ by $\varphi_n$ in \eqref{sect-4.6} and letting
$n\ra\infty$, we finally get the desired result. \fin

\end{document}